\documentclass{article}
\usepackage[normalem]{ulem}
\usepackage{amsthm}
\usepackage{thmtools}
\declaretheoremstyle[headfont=\normalfont]{normalhead}
\usepackage{color}
\usepackage{graphicx}
\usepackage{morefloats}
\usepackage{hyperref}
\usepackage{ amssymb }
\usepackage{textcomp}
\usepackage{url}
\usepackage{float}
\newtheoremstyle{mydef}
{\topsep}{\topsep}%
{}{}%
{\itshape}{}
{\newline}
{%
  \rule{\textwidth}{0.0pt}\\*%
  \thmname{#1}~\thmnumber{#2}\thmnote{\-\ #3}.\\*[-1.5ex]%
  \rule{\textwidth}{0.0pt}}%
\hypersetup{
   colorlinks=true,
    linkcolor=blue,
    filecolor=magenta,      
    citecolor=blue,
    urlcolor = blue
}
\begin{document}
\theoremstyle{mydef}
\newtheorem{conjecture}{Conjecture}
\newtheorem{theorem}{Theorem}
\newtheorem{question}{Question}
\newtheorem{remark}{Remark}
\newtheorem{proposal}{Proposal}
\newtheorem{lemma}{Lemma}
\newtheorem{corollary}{Corollary}

\author{Barry Brent}

\date{11h 11 August 2017}

\title{Variants of the Riemann zeta function}
\maketitle
\begin{abstract}
\hskip -.2in
We construct
variants of the Riemann zeta function 
with convenient properties
and make 
conjectures about their dynamics;
some of the conjectures are based on 
an analogy with the
dynamical system of zeta.
More specifically, we study the family
of functions
$V_z: s \mapsto  \zeta(s) \exp (zs)$.
We
observe convergence of 
$V_z$ fixed points
along nearly 
logarithmic spirals with
initial points at zeta 
fixed points and
centered upon Riemann zeros.
We can approximate these
spirals numerically,
so they
might afford a means to
study the geometry 
of the relationship
of zeta fixed points
to Riemann zeros.
\end{abstract}

\bibliographystyle{plain}

\section{\sc introduction}
In this article, we construct
variants of the Riemann zeta function 
with convenient properties
and make 
conjectures about their dynamics;
some of the conjectures are based on 
an analogy with the
dynamical system of zeta.
More specifically, we study the family
of functions
$V_z: s \mapsto  \zeta(s) \exp (zs)$.
We
observe convergence of 
$V_z$ fixed points
along nearly 
logarithmic spirals with
initial points at zeta 
fixed points and
centered upon Riemann zeros.
We can approximate these
spirals numerically,
so they
might afford a means to
study the geometry 
of the relationship
of zeta fixed points
to Riemann zeros.
\newline \newline
When we 
examined other zeta variants 
we observed behavior
similar to  that of the $V_z$;  we
do not have  a clear idea
of the extent of the phenomenon.
\newline \newline
We formulated  our conjectures
after  computer experiments using 
\it Mathematica. \rm Data and 
\it Mathematica \rm notebooks
for this project are at
ResearchGate, here \cite{Br3}.
Several other writers
have considered the dynamics of 
the Riemann zeta function,
 for example,
Kawahira \cite{Ka} and Woon \cite{W}.
\subsection{Definitions.}
In our experiments we found
evidence that  
certain sequences
of complex numbers are interpolated
by nearly logarithmic spirals,
but at first without independent
information about the 
underlying spirals.
In order to describe these 
observations,  we  use
terminology  that initially avoids
referring to  these
 spirals.
\newline \newline
Given a sequence of complex numbers
$\vec{z}$ 
$= \{z_0, z_1, ....\}$ 
with limit $\gamma$, we define a variant
$\theta_{\vec{z}}$
of the argument  concept
 by making
the following choices, which are 
always possible and
unique after fixing a branch of the 
argument function: \newline \newline
(1) $\theta_{\vec{z}}(z_0) = 
\arg(z_0-\gamma)$ \newline
(2) $\theta_{\vec{z}}(z_n) 
\equiv \arg(z_n - \gamma)$ mod $2\pi$, 
and
\newline
(3) for any non-negative integer $n$, 
$\theta_{\vec{z}}(z_{n+1})$ is 
the  smallest  real number
greater than  
$\theta_{\vec{z}}(z_n)$ 
compatible with
conditions (1) and (2). 
\newline \newline
Let $h$ and $K$ be positive integers. 
For each positive integer $n$ let 
$r(z_n) := |z_n-\gamma|$ and for 
$h = 0, 1, 2, ...$,
let $m, b$ be real numbers
such that the straight line
$y= m x + b$ is the best linear  model
of the data  
$(\theta_{\vec{z}}(z_n), \log r(z_n))$,
$n= h, h + 1, h + 2, ..., h+K$. 
(For the sake of 
definiteness,
we might insist that the model be 
chosen by the method of least squares; 
in practice, we have relied on
proprietary routines  of \it Mathematica\rm.)
\newline \newline
Now let  
$$d_h(n) :=  
\frac {|m  \theta_{\vec{z}}(z_n) + b - \log r(z_n)|}
{\log r(z_n)},$$
and let
$d_h := d_h(h+K)$. 
Then we  say that the sequence $\vec{z}$ 
is nearly logarithmic with respect to 
$\gamma$ if, for any choice of $K$,
$d_h  \to 0$ with exponential
decay as $h \to \infty$.
\newline \newline
Now we define  a nearly logarithmic spiral
in the complex plane 
as a spiral $S$ with center $\gamma$
such that any sequence of points
$\{z_n\}_n$ 
on $S$ 
with $|z_n - \gamma|$
decreasing monotonically with $n$
is nearly logarithmic with 
respect to $\gamma$.
\newline \newline
Next we say what it means for 
$\vec{z}$
 to
be  nearly uniformly distributed
with respect to $\gamma$; namely,
if $\delta_n$ is the quantity
$|\theta_{\vec{z}}(z_n) - \theta_{\vec{z}}(z_{n+1})|$
and  
$\Delta_n = |\delta_n - \delta_{n+1}|$,
then $\Delta_n$
decays exponentially with $n$.
\newline \newline
We  also adopt what is more or 
less standard usage, 
namely, that
an exactly logarithmic spiral with center  
$\gamma$ is a curve  consisting 
of points $z$ satisfying 
$\log r(z) =  m \theta_{\gamma}(z) + b$
for  real constants $m$ and $b$ (with the
above conventions for $\theta_{\gamma}$, 
except  that condition (3)
in the definition of  
$\theta_{\vec{z}}$ is replaced
with the requirement that $\theta_{\gamma}$ be 
continuous and monotonic
increasing as $z \to \gamma$.
An exactly logarithmic spiral
is (according to these definitions)
also a nearly logarithmic spiral.
In \cite{Br2}, 
we studied  the complex-valued
deviations of
certain nearly logarithmic
spirals  from  exactly logarithmic spirals
to which we had fitted them.  We  displayed
some plots of 
these deviations in 
Figure 6.1  of  \cite{Br2},
and we  make
a similar analysis in the discussion of 
conjecture 3
below.
\newline \newline
\textsection \thinspace 
We fix the following notation
for the remainder of the article: 
\newline \newline
$\bf{C}^{\it cut}$ := $\{z \in \bf{C}$ 
s.t. $z$
does not lie on the negative real axis$\}$.
\newline 
$\rho$ is a 
nontrivial Riemann zero.
\newline
$\Phi_z$ is the set of fixed points of $V_z$.
\newline $\psi$ 
is an element of  $\Phi_0$
(fixed point of zeta).
\newline $\psi_{\rho}$ is
the  particular
 element of $\Phi_0$ 
 closest to $\rho$.
\newline $u$ is
a point on the
unit circle.
\newline
$\vec{R}_u$ is the ray
emanating from zero and passing through $u$,
so that $z \in \vec{R}_u$
if and only if $z =xu$ for some
$x \geq 0$.
\newline
$X=(x_0, x_1,..., )$
is an infinite increasing
arithmetic  progression with 
$x_0 = 0$.
\newline 
$\vec{\phi} = (\phi_0, \phi_1, ...)$ 
is a member of
$\Phi_{0} \times \Phi_{x_1 u} 
\times \Phi_{x_2 u} \times ... := \vec{\Phi}_{X,u}$.
\newline
$V^g_z(s):= g(s)e^{zs}$, 
$\Phi^g_z$ is the set of
fixed points of $V^g_z$,
and
$\Phi^g_{0} \times \Phi^g_{x_1 u} 
\times \Phi^g_{x_2 u} \times ... := \vec{\Phi^g}_{X,u}$.
(Thus $V_z = V^{\zeta}_z$, etc.)
\subsection{A theorem on limits 
of sequences from
\texorpdfstring{$\Phi_{X, u}$}
{TEXT}.} 
Before we list the conjectures,
we prove a theorem.
\begin{theorem} 
Suppose that
$\vec{\phi}$ 
in $\Phi^g_{X, u}$
converges to a complex
number $\lambda$,
$g$ is continuous
at $\lambda$,
and that
$\Re(\lambda) \cdot \Re(u)
-\Im(\lambda) \cdot \Im(u) > 0$. 
Then $g(\lambda) = 0$.
\end{theorem}
\hskip -.2in  
\it Proof  \rm 
\thinspace
Let us write $\lambda = L + Mi$
and $u = P + Qi \thinspace (L, M, P, Q$ real.)
If $g(\lambda) \neq 0$, then
we can choose a subsequence 
$\vec{\phi^*}$ of $\vec{\phi}$
and a positive number $B$ 
such that for, some
natural number $n(B)$,
$n > n(B) \Rightarrow
|g(\phi^*_n)| > B$.
Let us write the real
and imaginary parts of
$\phi^*_n$ as $a_n, b_n$, respectively,
and set 
$D_n = a_n P - b_n Q, E_n = a_n Q + b_n P$.
From the hypotheses we know that
the $D_n \rightarrow LP - MQ > 0$,
and so there is a 
number $C > 0$ and a number 
$n(C)$ such that
$n > n(C) \Rightarrow D_n > C > 0 
\Rightarrow e^{D_n} > e^C > 1$.
Let us choose a subsequence 
$X^* = (x^*_0, x^*_1, ...)$ of 
the arithmetic progression $X$
so that the $\phi^*_n$ are fixed points
of the functions $V^g_{x^*_n u}$.
Thus $\phi^*_n = V^g_{x^*_n u} (\phi^*_n)= 
g(\phi^*_n) e^{x^*_n u \phi^*_n}$.
Now we have:
$$|e^{x^*_n u \phi^*_n}| =
|e^{x^*_n (P + Qi) (a_n + b_n i)}|
= |e^{x^*_n  D_n + i x^*_n E_n} | = e^{x^*_n D_n}.$$
Therefore,
$$|\phi^*_n| = |g(\phi^*_n)| \cdot 
|e^{x^*_n u \phi^*_n}| = |g(\phi^*_n)|e^{x^*_n D_n}.$$
Now let $n > \max (n(B), n(C))$; it follows that
$|\phi^*_n| > B (e^C)^{x^*_n}$. 
But $x^*_n \rightarrow \infty$,
hence,
$(e^C)^{x^*_n} \rightarrow \infty$.
Since the $\phi^*_n$ converge to the finite number 
$\lambda$, this cannot be true.
$\square$ \newpage
\begin{corollary} 
Suppose that
$\vec{\phi}$ 
in $\Phi_{X, u}$
converges to a complex
number $\lambda \neq 1$
and that \newline
$\Re(\lambda) \cdot \Re(u)
-\Im(\lambda) \cdot \Im(u) > 0$. 
Then $\lambda$ is a Riemann zero.
\end{corollary}
\subsection{Conjectures.}
In all of the conjectures to follow,
we assume that $u \neq -1$ and that 
(if $u \neq 1$) 
$\Im \rho \cdot \Im u < 0$. 
Our experiments indicate that
this restriction is necessary.
\begin{conjecture} \rm 
\vskip .1in
(1) The  set of imaginary parts 
$\{\Im \phi: \phi \in \Phi_{z} \}$ 
is unbounded
and nonempty.
\newline 
(2) A unique $\psi_{\rho}$
exists for each $\rho$.
\end{conjecture}
\begin{conjecture} \rm
(1) Some  $\vec{\phi}$ in
$\vec{\Phi}_{X,u}$
converges 
to each $\rho$. 
(2) Any such $\vec{\phi}$ 
is nearly logarithmic 
with respect to  $\rho$,
 and nearly uniformly distributed 
with respect to $\rho$.
\end{conjecture}
\begin{conjecture} \rm
For each choice
of $\rho$,
there is a continuous function
$f^{cut}: \bf{C}^{\it cut} 
\rightarrow \bf{C}$
such that the restriction of 
$f^{cut}$ to $\vec{R}_u$
is a function $f$ with the
following properties:
\newline
(1) 
$f: \vec{R}_u \rightarrow \bf{C}$ 
is continuous and one-to-one,
\newline
(2) $f(z) \in \Phi_{z}$
for each $z \in \vec{R}_u$,
\newline
(3)  $f(0) = \psi_{\rho}$,
\newline
(4)  $\lim_{z \to \infty} 
f(z) =  \rho$,
\newline
(5) the image of $\vec{R}_u$ 
under $f$ is a nearly
logarithmic
spiral $S_{u,\rho}$ with center 
$\rho$,
\newline
(6)  among the sequences
$\vec{\phi} = 
(\phi_0, \phi_1, ...) \in 
\vec{\Phi}_{X,u}$
converging to $\rho$, 
all the points of one of them
(say, $\vec{\phi}_{u,\rho}$)
are interpolated by $S_{u,\rho}$;
the initial element of  
$\vec{\phi}_{u,\rho}$
is  
$\psi_{\rho}$.
\newline
(\it n.b. \rm Our notation suppresses the dependence of
$f$ on $u$ and $\rho$, and
the dependence of $f^{cut}$
on $\rho$.)
\end{conjecture}
\hskip -.2in
\textsection  \thinspace 
The next conjecture is
motivated by the analogy mentioned in the 
introduction and discussed later in the article.
\begin{conjecture} \rm
For each $X$ and $u$,
there is a function $G$ with the following
properties:\newline
(1) $G$ is meromorphic
and independent of $\rho$.
\newline
(2) The sequence  $\vec{\phi}_{u,\rho}$
from conjecture 3
satisfies
$G(\phi_n) =  \phi_{n-1}$
for  $n = 1, 2, ....$ 
\newline
(3)  Each $\rho$ is a repelling fixed point
of $G$,
\newline
(4) $G$ is many-to-one,
but for each $\rho$  there 
is a function 
$F_{\rho}$
such that 

(i)  $F_{\rho}^{(-1)} = G$,

(ii) consequently (in view of clause 2) 
for $n = 0, 1, ...$,
$$F_{\rho}(\phi_n) =  \phi_{n+1},$$

(iii) 
 $\rho$ is
an attracting fixed point of 
$F_{\rho}$,
and 

(iv) $F_{\rho}$ 
carries $S_{u,\rho}$
into itself.
\newline
(\it n.b. \rm Our notation suppresses the 
dependence of $G$,
and the dependence of
$F_{\rho}$ on $X$ and $u$. 
\end{conjecture}
\begin{remark}
Obviously, it is 
 not surprising 
 that the $V_z$ might in some
way pick out
zeta zeros and fixed points
as special.
On the other hand, 
our numerical methods make it
feasible
to estimate the 
parameters incorporated in
the equations of
simple curves interpolating 
$V_z$ fixed points. 
This is possible because we are 
able to find points on the
curves that are arbitrarily 
close to each other, something
we were unable to do for
analogous curves we studied in \cite{Br2},
where to do so would have required 
doing something else we do not know
how to do: extending
the iteration of zeta to non-integer heights.
(In the notation we will introduce below,
it would have required evaluating
expressions of the form
$\zeta^{\circ q}(s)$
for non-integer values of $q$.
There is some literature around the
problem of extending the
iteration operator, \it e.g.\rm, 
\cite{N2, N1}, but we have not
reduced those
results to code. In the
situation of the present article,
we do not have to
confront this difficulty.)
In the present situation,
conjecturally, 
these curves are nearly logarithmic 
spirals with 
a Riemann zero $\rho$ 
and the zeta fixed point 
$\psi_{\rho}$ 
as center and initial point, 
respectively. At present,
we are able to make
the estimates only for
one $\rho$
at a time. But
more study of these
estimates may be a way
to search for a
dictionary between
zeta zeros and zeta fixed points.
The Riemann hypothesis might
then be rephrased as a 
claim about the fixed
points. In another direction,
the convergence properties
of sequences from the
$\Phi_z$
bear upon the 
Riemann hypothesis
as well. (See question 1 below.)
\end{remark}
\begin{remark}
Conjecture 3 
(in which $\psi = \psi_{\rho}$) 
is
supported by some
experimental evidence
(see below). Based on the
analogy with the
situation of \cite{Br2} 
promised in the 
introduction and discussed in
the next section,
it seems plausible 
that conjecture 3
might extend to
arbitrary $\psi$.
To find experimental evidence
for other $\psi$
(corresponding to
evidence for the
analogous claim in \cite{Br2}),
we  would first need to identify 
the unknown
zeta-analogue $G$, because
the needed sequences in \cite{Br2} 
corresponding to the
$\vec{\phi} \in \vec{\Phi}_{X,u}$
in  conjectures  3 and 4
were obtained by  solving equations 
involving zeta.
Conjecture 4 is also
a claim about $\psi_{\rho}$,
since it involves the
curve $S_{u, \rho}$ and 
clause (6) of conjecture 3
associates $\psi_{\rho}$
to $S_{u, \rho}$.
Therefore the question arises
of the extension
of conjecture 4 to
arbitrary $\psi$.
Conjecture 4 is
based entirely on the analogy
with the situation of \cite{Br2}.
Thus, the same requirements
for finding experimental
evidence for its extension
to arbitrary $\psi$ apply
to conjecture 4 itself.
\newline \newline
We write down the
contemplated extensions 
of
conjectures 3 and 4
(labeled more tentatively there as
``proposals'')
in more detail in the 
appendix.
\end{remark}
\vskip .05in
\hskip -.2in
\textsection  \thinspace 
When iteration of
a function $f$ is meaningful,
we write
$f^{\circ 0}(s) = s$ and, for  
$n$ a positive integer,
$f^{\circ n}(s) = f(f^{\circ (n -1)}(s))$.
\begin{remark} 
For particular $\rho$,
by conjecture 1 we
may choose  
$\psi_{\rho} \in \Phi_0$.
Then conjecture 4
 implies that 
the  forward orbit of $\psi_{\rho}$ 
under $F_{\rho}$, namely
\newline
$(\psi_{\rho},  F_{\rho}(\psi_{\rho}),
F_{\rho}^{\circ 2}(\psi_{\rho}), ...)$,
is  a sequence
in  $\vec{\Phi}_{X,u}$
converging to  $\rho$.  Therefore,
clause (1) of
conjecture 2 follows from  
conjecture 1 and conjecture 4.
\end{remark}
\begin{conjecture} \rm
Here, it will be
convenient to indicate the
dependence of $G$ upon
$X$ and $u$ 
by writing $G = G_{X, u}$.
Suppose that
$X'$ is a refinement of  
$X$ in the sense that
the common difference  
between consecutive elements of  $X$ 
and between  consecutive elements of 
$X'$ are  $d, d'$ respectively
and  $d=Kd'$ for some
natural number $K$. Then 
$G_{X,u} = G_{X',u}^{\circ K}$.
\end{conjecture}
\subsection{Conjectures 4, 5
and the promised analogy.} 
In \cite{Br2}, we studied
the dynamical system of zeta
and several associated 
nearly logarithmic 
sequences
$\vec{\phi}_{\psi, \rho, \zeta} = 
(z_0, z_1, ...)$
such that \newline \newline
(*) $z_0 = \rho$, \newline
(**) $z_{n-1} = \zeta(z_n)$ for all 
$n > 0$,
and 
\newline
(***) 
$\lim \vec{\phi}_{\psi, \rho, \zeta} = \psi$.
\newline \newline
In that article, we
conjectured that such
a sequence existed
for each choice of the 
pair $(\psi, \rho)$,
and we remarked that
in view of (for example) 
\cite{HY}, Theorem 2.6,
this phenomenon
would follow if
each $\psi$ is a
repelling fixed point
of zeta, hence an attracting
fixed point of a branch of the 
inverse $\zeta^{(-1)}$.
(Here \it branch \rm has
its usual meaning in 
complex analysis;
it does not denote a branch of a 
backward orbit, as defined in 
\cite{Br2},
but it is true that  
$\vec{\phi}_{\psi, \rho, \zeta}$
constitutes a branch of the 
backward orbit of zeta as
defined there.)
The behavior of iterates
of zeta were shown to
be similar in the range of our observations.
So we think it
is plausible to
propose a dynamical explanation
for the attraction of
the sequences in $\vec{\Phi}_X$
to Riemann zeros that draws on an analogy
with the proposals in \cite{Br2}.
In this analogy, in the 
context of \cite{Br2},
with 
 $X_n$ (say) $= (0, n, 2n, ...)$
$G$ should
correspond to $\zeta^{\circ n}$, 
the $F$'s
should
correspond to (complex analysis sense) 
branches of $\zeta^{\circ n (-1)}$,
and
the roles of the 
zeta fixed points $\psi$ and the
zeta zeros $\rho$ in \cite{Br2}
should be swapped.
This idea results in conjecture 4,
conjecture 5,
and proposal 2 (appendix.)
\newline \newline
We ask
the following question,
because (together with Conjecture 2)
an affirmative answer implies
the truth of the Riemann hypothesis,
and (together with Theorem 1)
a negative answer provides an avenue
to search for counterexamples
to the Riemann hypothesis.
Both possibilities are consistent with
our observations.
\begin{question}
Is the following claim true?
 \it Suppose $0 < \sigma < 1$
 and the real parts 
$(\Re \phi_0, \Re \phi_1, ...)$
converge to  $\sigma$. Then
$\sigma = \frac 12$.
\end{question}
\section{\sc methods}
\subsection{Quadrant plots.}
\begin{figure}[!htbp]
\centering
\includegraphics[scale=.8]{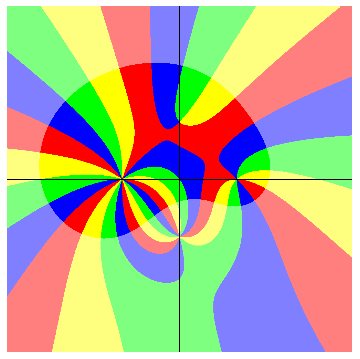}
\vskip .07in
\caption{Quadrant plot of $s \mapsto 
(s - 1)^2 (s - i) (s + 1)^5/(s+i)^3$}
\label{leaf2point1.png}
\end{figure}

We will make use 
of colored plots
(``quadrant plots'').
The  routine that makes
quadrant plots
takes as inputs
a specification 
of the image resolution,
the center and dimensions
of a region $R$
in the complex plane,
and a routine to compute some
$\bf{C} \rightarrow \widehat{\bf{C}}$ 
function
$f$. The
output colors a 
small square (the size of
which depends on the 
resolution) 
around each point $w$
of a regular lattice in $R$;
the color is chosen to
represent
the quadrant of $f(w)$.  
Similar methods
that do not use coloring are
put to use in, \it e.g\rm,
\cite{A}.)
\newline \newline
The pixel representing the square
is colored according to the rules in 
Table 1.
In the table,
the region $D$ is a disk with center $s = 0$ and 
large radius $r$ (chosen as may be convenient.) 
We denote the complement of $D$ as $-D$.
\newline \newline 
\hskip 1in
\begin{tabular}{|c|c|} \hline
Location of $f(s)$& Color of pixel depicting region containing 
$s$\\ \hline \hline
real and imaginary axes & black \\ \hline 
$D \hskip .05in\cap $ Quadrant I  &  rich blue\\ \hline 
$ - D \hskip .05in\cap$ Quadrant I  &  pale blue\\ \hline
$D \hskip .05in\cap $ Quadrant II  & rich red \\  \hline
$ - D \hskip .05in\cap$ Quadrant II  &  pale red\\ \hline
$D \hskip .05in\cap $ Quadrant III  & rich yellow \\ \hline 
$ - D \hskip .05in\cap$ Quadrant III  &  pale yellow\\ \hline
$D \hskip .05in\cap $ Quadrant IV  &  rich green\\ \hline
$ - D \hskip .05in\cap$ Quadrant IV  &  pale green\\ \hline 
\end{tabular} 
\vskip .1in
\hskip .5in
\sc{table 1: coloring scheme for quadrant plots}
\rm
\vskip .1in 
\hskip -.2in
The junction of four rich colors 
represents a zero,
the junction of four pale colors represents a pole,
and
the boundary of two appropriately-colored regions is an 
$f$ pre-image of an
axis. An example is shown in Figure \ref{leaf2point1.png}.
(We have superimposed a pair of axes on this 
quadrant plot.)
\newpage
\subsection{Basins of attraction.}
For $z \in \bf{C} \rm \cup \{\infty\}$, 
let $A_z :=\{w \in \bf{C}$ s.t.
$\lim_{n \to\infty} 
\zeta^{\circ n}(w) = z\}$
(the ``basin of attraction'' of $z$ under 
zeta iteration.)
Then $A_{\infty}$ and its
complement in $\bf{C}$ are fractals \cite{W}.
(As far as we know,
plots of  $A_{\infty}$
were made first by 
Woon in \cite{W}.)
Let $\phi \approx -.295905$
be the largest negative
zeta fixed point.
As we noted in \cite{Br2},
plots of
$A_{\phi}$ and the complement in
$\bf{C}$ of $A_{\infty}$
are indistinguishable to the eye.
(But, for example,
the complement in
$\bf{C}$ of $A_{\infty}$
contains the zeta cycles and zeta zeros, and
$A_{\phi}$ does not.)
We reproduce plots of 
$A_{\phi}$ and the complement in
$\bf{C}$ of $A_{\infty}$
in Figure \ref{leaf1point1.png}
(appendix.)
\section{\sc Basis of the claims}
\subsection{Conjecture 1.}
\subsubsection{Existence and uniqueness of 
\texorpdfstring{$\psi_{\rho}$}{TEXT}.}
The right panel of Figure 
\ref{leaf3point2left.png}
(Figure 3.2 of \cite{Br2}) 
in the appendix is 
a quadrant plot
of the function 
$s \mapsto \zeta(s) - s$,
the zeros of which are the zeta fixed points.
The quadrant plot is shown superposed on
a plot of $A_{\phi}$.
The major feature
of this escaping set 
is that it consists of an infinite family of
irregularly shaped bulbs straddling
the critical strip, each
one of which 
(except for the central bulb, referred to
in \cite{Br2} as a
``cardioid'') contains one Riemann
zero (trivial or nontrivial,
but only the cardioid and 
the bulbs associated to
nontrivial zeros are visible at the scale
of Figure \ref{leaf3point2left.png}.) 
One zero of 
$s \mapsto \zeta(s) - s$
evidently lies on
one filament decorating
each of
the visible 
non-cardioid bulbs. 
These filaments
each consists of smaller bulbs
of the zeta escaping set, and
there is a numerical pattern
(described in \cite{Br2})
dictating the distribution of zeta 
fixed points among 
these bulbs. The fact that
this distribution is non-random
is evidence (we would argue)
that the association of
zeta fixed points to the Riemann zeros
is one-to-one, and, therefore,
that  there is probably
a unique $\psi_{\rho}$
for each $\rho$.
\subsubsection{The case
\texorpdfstring{$u=1$}
{TEXT}.} 
Clause 1 of conjecture 1 is the claim that  
the  set of imaginary parts 
$\{\Im \phi: \phi \in \Phi_z \}$ 
is unbounded
and nonempty.
The
set $\Phi_z$ is the set of solutions 
of the function 
$s \mapsto V_z(s) - s$.
When $z = 0$, 
$\Phi_z = \Phi_0$ is the set of fixed points
of the Riemann zeta function. 
Plots of these fixed points in \cite{Br2} 
provided our evidence for this case of 
the conjecture; we will not reproduce
them here. 
Using our particular
methods, we can only spot-check
this claim for selected $z$
and selected regions of the complex
plane. 
In Figure \ref{seqfig2.png} below, we display a
quadrant plot of the function 
$s \mapsto V_1(s) - s$
on a $2000$ by $2000$ square
with center at $s=0$ in the complex plane.
\begin{figure}[!htbp]
\centering
\includegraphics[scale=.8]{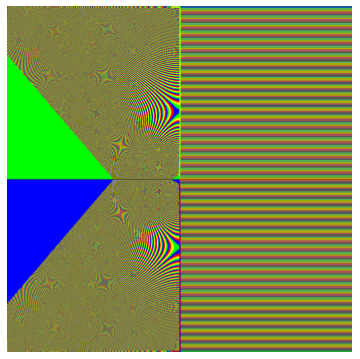}
\vskip .07in
\caption{2000 by 2000 quadrant plot of $s 
\mapsto V_1(s) - s$}
\label{seqfig2.png}
\end{figure}
Visible zeros of this function (lying at 
the junctions of four colors) lie in the vicinity
of the imaginary axis. There are also zeros (not visible
at this scale) in the vicinity of each nontrivial
Riemann zero in the depicted region,
on the boundary of a triangular region 
(barely
visible at this scale) near the center of the plot, 
and along the 
negative real axis.
The straight
boundaries of the green and blue
regions  do 
not (on inspection at smaller scales)
 contain any zeros of this function.
 \newline \newline
 In Figure \ref{seqfig3r1c1.png} below, we  display
quadrant plots of 
$s \mapsto V_k(s) - s, k = 2, 4, 8$
and  $16$, each centered at $s = 0$,
in squares with side length  
$400,  2400, 1.2 \times 10^5$, and  
$2 \times 10^8$ respectively.
The are shown clockwise from upper left in the 
order of the size of $k$.
 \begin{figure}[!htbp]
\centering
\includegraphics[scale=.45]{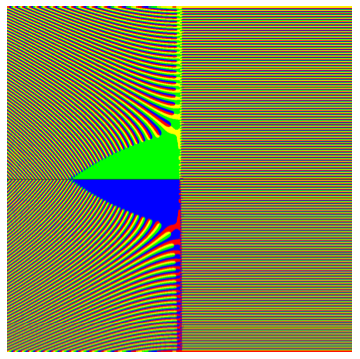}
\hskip .05in
\includegraphics[scale=.45]{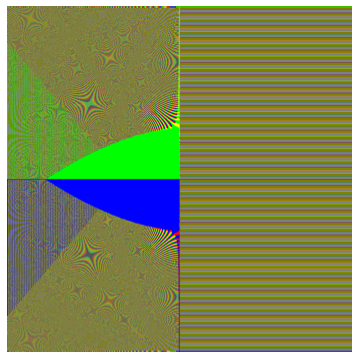}
\vskip .05in
\includegraphics[scale=.45]{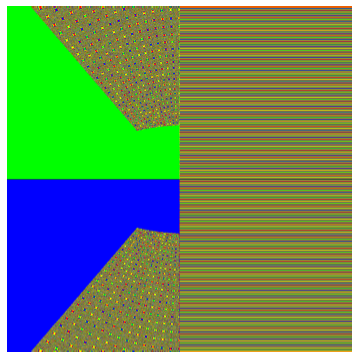}
\hskip .05in
\includegraphics[scale=.45]{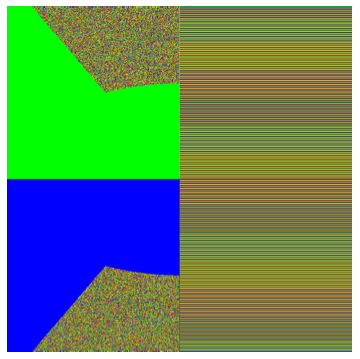}
\vskip .07in
\caption{Quadrant plots of $s 
\mapsto V_k(s) - s, k = 2, 4, 8$ and $16$}
\label{seqfig3r1c1.png}
\end{figure}
\newline \newline
Since the Riemann hypothesis has been
verified well beyond
the range of our experiments,
we are safe in designating the
$n^{th}$-by-height nontrivial zero  
in the upper half plane as $\rho_n$.
Figure \ref{seqfig4left.png} 
shows (at the four-color junctions)
zeros of 
$s \mapsto V_{100}(s) - s$ on two 
$1.2$ by $1.2$ squares, 
one centered at $\rho_1$
and the other centered at $\rho_{649}$.
The zeros along the the 
left sides of the squares
apparently belong (as we will explain below) 
to sequences of fixed points  $\phi_n$
of the $V_n$, the real parts of 
which converge to zero
(but it is not clear that 
the imaginary parts converge at all.)
The zeros very near 
the centers of the squares
appear to belong
to such sequences converging to $\rho_1$
and $\rho_{649}$, respectively,
as we will also explain below.
 \begin{figure}[!htbp]
\centering
\includegraphics[scale=.4]{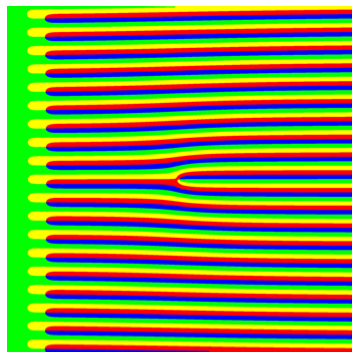}
\includegraphics[scale=.4]{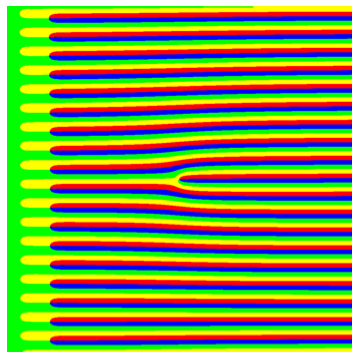}
\vskip .07in
\caption{Quadrant plots of $s \mapsto V_{100}(s) - s$
near $\rho_1$ and $\rho_{649}$.}
\label{seqfig4left.png}
\end{figure}
These and similar plots, which depict the 
fixed points of various $V_z$ as (apparently) 
unbounded sets of isolated points
in the complex plane,
are the basis of conjecture 1.
\subsection{Conjecture 2: the case
\texorpdfstring{$u=1$}{}.}
This is the
claim that, for $X$ and $\rho$
as above,
some  $\vec{\phi} \in \vec{\Phi}_{X,1}$
converges to $\rho$,
is nearly logarithmic with respect to
$\rho$, and is nearly uniformly 
distributed  with  respect to $\rho$.
We will describe experiments 
that tend to support this claim
for
 $X = (0, 1, 2, ...), X' = 
 (0, \frac 15, \frac 25, ...)$ 
 and for  
$\rho_n$ with $n$ selected from 
the range $1 \leq n \leq 100$. 
\subsubsection{Typical plots.}
\begin{figure}[!htbp]
\centering
\includegraphics[scale=.45]{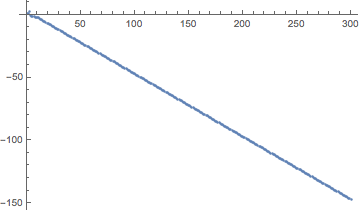}
\hskip .05in
\includegraphics[scale=.45]{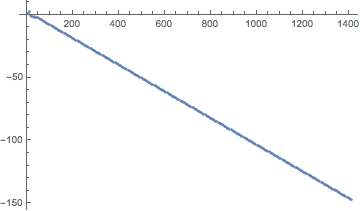}
\vskip .05in
\includegraphics[scale=.45]{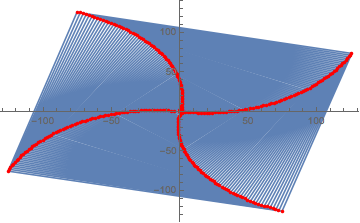}
\hskip .05in
\includegraphics[scale=.45]{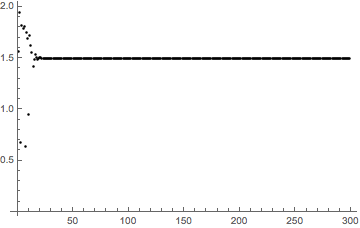}
\vskip .05in
\includegraphics[scale=.45]{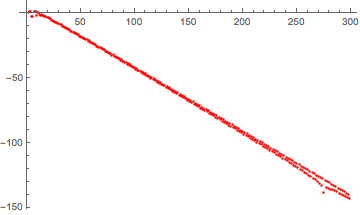}
\hskip .05in
\includegraphics[scale=.45]{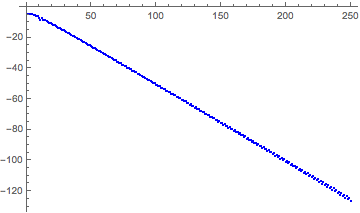}
\vskip .07in
\caption{Zero $\rho_1, u = 1, X  = (0, 1, 2, ...)$, 
$0 \leq n \leq 300$
\newline row 1, column 1:
$\log |\phi_n - \rho_1|$ vs. $n$
\newline  row 1, column 2: 
$\log |\phi_n - \rho_1|$ vs. $\theta_{\vec{\phi}}(\phi_n)$
\newline row 2, column 1: logarithmically scaled plot
of  $\vec{\phi}$
\newline row 2, column 2: $\delta_n/\pi$ vs. $n$
\newline 
\newline row 3, column 1: $\log \Delta_n$ vs. $n$
\newline row 3, column 2:  $\log h_n$ vs. $n$ (K =  50)}
\label{seqfig5r1c1.png}
\end{figure}

\begin{figure}[!htbp]
\centering
\includegraphics[scale=.45]{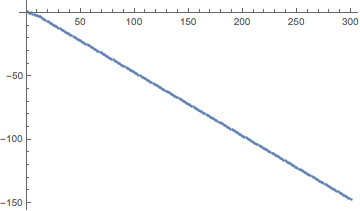}
\hskip .05in
\includegraphics[scale=.45]{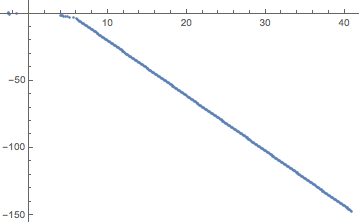}
\vskip .05in
\includegraphics[scale=.45]{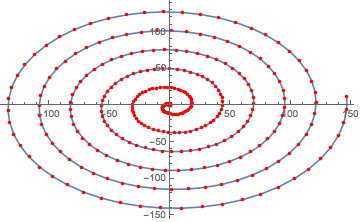}
\hskip .05in
\includegraphics[scale=.45]{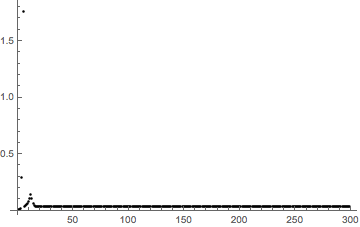}
\vskip .05in
\includegraphics[scale=.45]{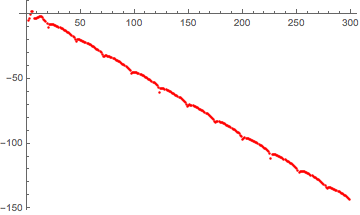}
\hskip .05in
\includegraphics[scale=.45]{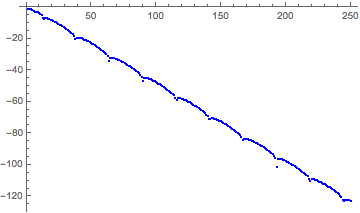}
\vskip .07in
\caption{Zero $\rho_3, u = 1, X  = (0, 1, 2, ...)$, 
$0 \leq n \leq 300$
\newline row 1, column 1:
$\log |\phi_n - \rho_1|$ vs. $n$
\newline  row 1, column 2: 
$\log |\phi_n - \rho_1|$ vs. $\theta_{\vec{\phi}}(\phi_n)$
\newline row 2, column 1: logarithmically scaled plot
of  $\vec{\phi}$
\newline row 2, column 2: $\delta_n/\pi$ vs. $n$
\newline 
\newline row 3, column 1: $\log \Delta_n$ vs. $n$
\newline row 3, column 2:  $\log h_n$ vs. $n$ (K =  50)}
\label{seqfig6r1c1.png}
\end{figure}

The  red points
on the pictured spirals
depict elements of the 
sequences $\vec{\phi}$ from 
conjecture 3, clause 5. 
Because the $\vec{\phi}$
converge to $\rho$ so rapidly,
logarithmic scaling was
necessary  to obtain readable plots.
As a result, the spirals are 
for the most part everted:
points that appear farther from the center
(which means farthest from the zero)
are, in fact, closer to the center.
(The only exception might be  the
red point  depicting $\psi_{\rho}$,
since  typically 
$|\psi_{\rho} - \rho| > 1$.)
Thus, 
for $\phi \in \vec{\phi}$,
the corresponding red point 
$p$ (say)
in the plots of spirals
below  is situated
on  the  ray  originating
at $\rho$ and passing through
$\phi$, but 
$|p -  \rho|  = \log |\phi - \rho|$.
In
these figures, blue chords
connect  representatives (red points)
of consecutive members of 
$\vec{\phi}$. 
\newline \newline
We  made plots  from which we  formulated 
conjecture 2. In one experiment,
we  
examined the  zeros 
$\rho_n, 1 \leq n \leq 100$.
Figures \ref{seqfig5r1c1.png} and 
\ref{seqfig6r1c1.png}  are typical outcomes
for $X = (0, 1, 2, ...)$. 
We chose them for display because,
together, they  indicate the dependence
upon $\delta_n$ of the  appearance
of the spirals.  
Figure  7  depicts several 
spirals centered on other zeros; 
we omit the  statistics
for these zeros, which are 
consistent with our conjectures.
\newpage
 \begin{figure}[!htbp]
\centering
\includegraphics[scale=.40]{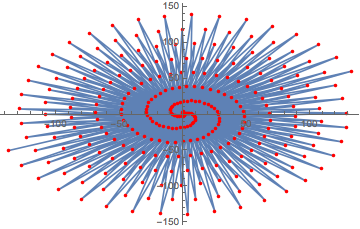}
\hskip .03in
\includegraphics[scale=.40]{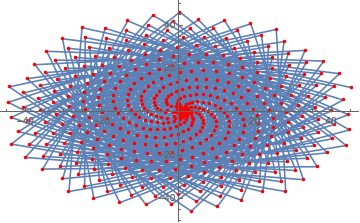}
\vskip .03in
\includegraphics[scale=.40]{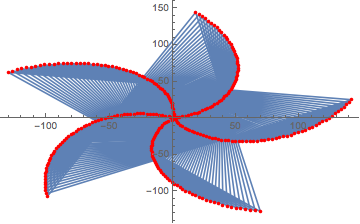}
\hskip .03in
\includegraphics[scale=.40]{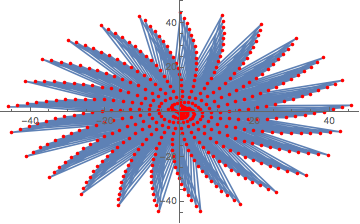}
\vskip .03in
\includegraphics[scale=.40]{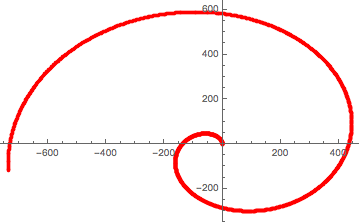}
\hskip .03in
\includegraphics[scale=.40]{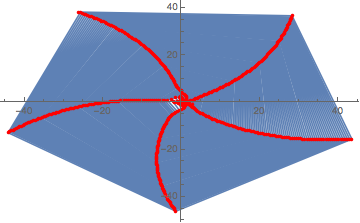}
\vskip .03in
\includegraphics[scale=.40]{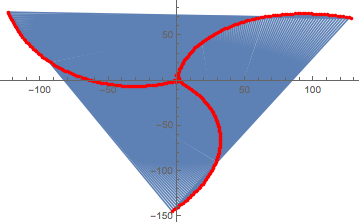}
\hskip .03in
\includegraphics[scale=.40]{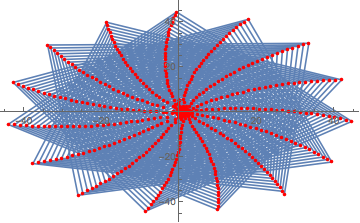}
\vskip .03in
\includegraphics[scale=.40]{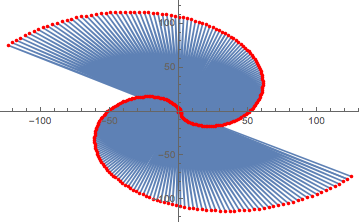}
\hskip .03in
\includegraphics[scale=.40]{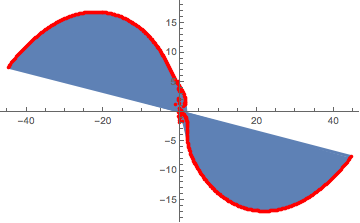}
\vskip .03in
\caption{Row 1, $\rho_7$; row 2, $\rho_{65}$;
row 3,  $\rho_{70}$; row 4, $\rho_{78}$,
row 5, $\rho_{82}$;\newline
column 1, $X = (0, 1, 2, ...)$; column 2, 
$X = (0, \frac 15, \frac 25, ...)$}
\label{seqfig7zr7left.png}
\end{figure}
\begin{figure}[!htbp]
\centering
\includegraphics[scale=.45]{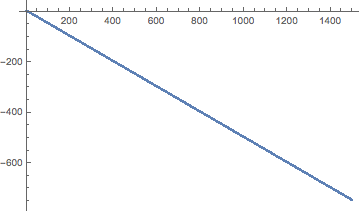}
\hskip .05in
\includegraphics[scale=.45]{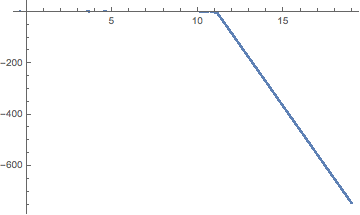}
\vskip .05in
\includegraphics[scale=.45]{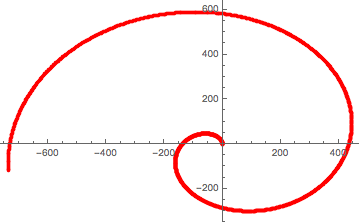}
\hskip .05in
\includegraphics[scale=.45]{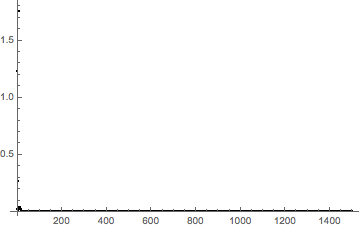}
\vskip .05in
\includegraphics[scale=.45]{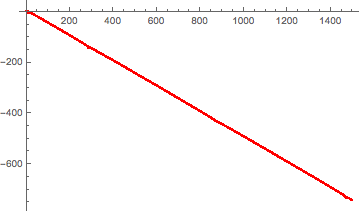}
\hskip .05in
\includegraphics[scale=.45]{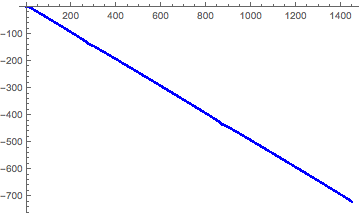}
\vskip .07in
\caption{Nontrivial zero $\rho_{70}, X  = (0, 1, 2, ...)$, 
$0 \leq n \leq 1500$;\newline
row 1, column 1:
$\log |\phi_n - \rho_1|$ vs. $n$;\newline
row 1, column 2: 
$\log |\phi_n - \rho_1|$ vs. $\theta_{\vec{\phi}}(\phi_n)$;\newline
row 2, column 1: logarithmically scaled plot
of  $\vec{\phi}$;\newline
row 2, column 2: $\delta_n/\pi$ vs. $n$;\newline
row 3, column 1: $\log \Delta_n$ vs. $n$;\newline
row 3, column 2:  $\log h_n$ vs. $n$ (K =  50)}
\label{seqfig8r1c1.png}
\end{figure} \hskip -.1in
\subsubsection{An anomaly.}
For the zero
$\rho_{70}$, while the other
statistics we display in these 
plots were entirely consistent
with our conjectures, $\delta_n$
was so small in the range 
$0 \leq n \leq 300$ that the spiral 
structure for $X = (0, 1, 2, ...)$ 
was not obvious and
we had to extend our observations 
to the range $0 \leq n \leq 1500$
to see it. 
It appears that,  for 
$\rho_{70}$,  $\delta_n$
converges rapidly to a limit $\approx .0016 \pi$,
so that  it requires  $1187$ points of 
$\vec{\phi}$ to wind entirely around
the  zero.
We have not encountered 
another spiral like this. 
(The anomaly does not extend, for example,
to the
case $X = (0, \frac 15, \frac 25,...)$.)
The uniqueness (within the 
range of our observations) 
of this anomaly suggests that,
for the purpose of searching
for the unknown functions $G$,
one should look for an invariant of
the Riemann zeros that takes
an anomalous value at $\rho_{70}$.
As we said above,
we regard
the Riemann zeta function 
as analogous to $G$;
spirals  studied in 
\cite{Br2} are centered
upon zeta fixed points $\psi$,
and  $\lim_{n \to \infty} \delta_n$ 
appears to be determined by
$\arg \frac d{ds}\zeta(s)|_{s = \psi}$.
Therefore 
 we might expect
 the  invariant we are looking for
to
coincide with 
 $\arg \frac d{ds}G(s)|_{s = \rho_n}$,
which we might hope to identify 
without knowing $G$ explicitly.
We  have undertaken a 
cursory search for 
such an
invariant, so far without success.
The plots for
$\rho_{70}$ are in Figure \ref{seqfig8r1c1.png}. 
An interesting feature of the spiral
about $\rho_{70}$ for
$X = (0, 1, 2, ...)$ is that 
red points  $p_n$ and $p_{n'}$ on it
are close just if  $|n - n'|$
is small; on most  of the other
spirals we will display
in this article, this is not the case,
as one can  tell  by  keeping track of the
blue connecting chords. 
(The spiral about $\rho_3$
depicted in Figure \ref{seqfig6r1c1.png},
for which the $\delta_n$ are
also fairly small, is an 
exception.) In the plot
shown in row 3, column 1 of 
Figure \ref{seqfig7zr7left.png} 
(and reproduced in Figure 
\ref{seqfig8r1c1.png}),  
consecutive  $p_n$ are
so close that  these chords are
not visible.
\subsection{Conjecture 3.}
\begin{figure}[!htbp]
\centering
\includegraphics[scale=.45]{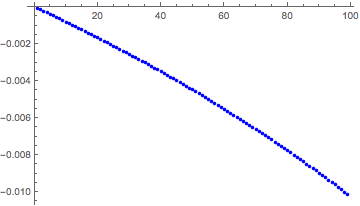}
\hskip .05in
\includegraphics[scale=.45]{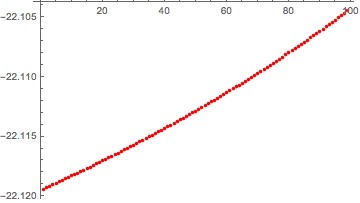}
\vskip .05in
\includegraphics[scale=.45]{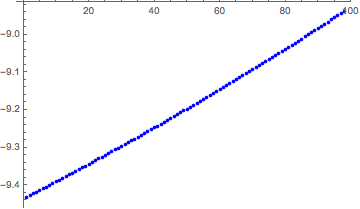}
\hskip .05in
\includegraphics[scale=.45]{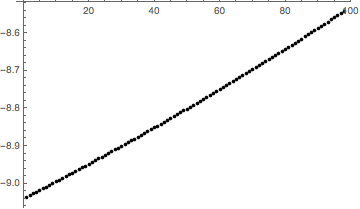}
\caption{Linear models of $\log r$ vs.
$\theta$ for $\rho = \rho_1, u = 1$,
and a sequence
$\vec{\phi} \in \Phi_X$ on coarsened, truncated
$X_1 = (50, 50.001, 50.002, ..., 100)$: parameters plotted against
the filter index.
\newline
row 1, column 1: slopes;
\newline 
row 1, column 2: intercepts;
\newline
row 2, column 1: logarithms of slope increments;
\newline
row 2, column 2: logarithms of intercept increments
}
\label{seqfig9r1c1.png}
\end{figure}

\begin{figure}[!htbp]
\centering
\includegraphics[scale=.45]{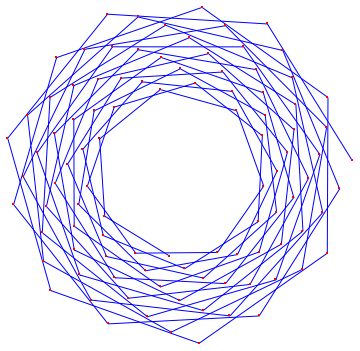}
\hskip .05in
\includegraphics[scale=.45]{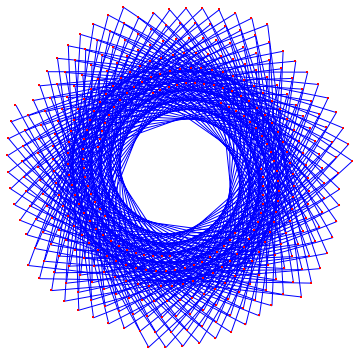}
\vskip .05in
\includegraphics[scale=.45]{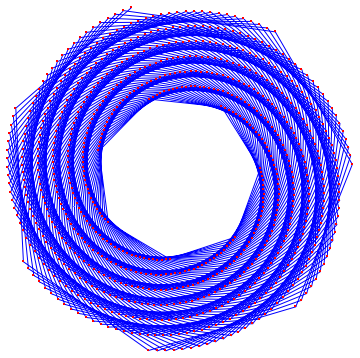}
\hskip .05in
\includegraphics[scale=.45]{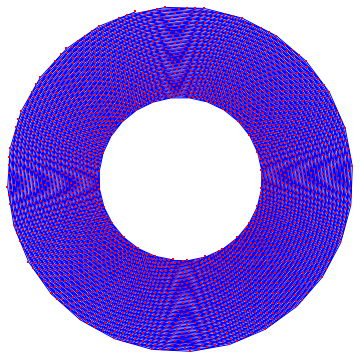}
\vskip .05in
\includegraphics[scale=.45]{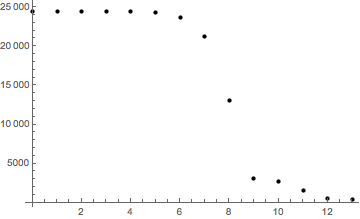}
\hskip .05in
\includegraphics[scale=.45]{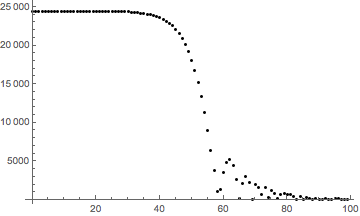}
\caption{$\rho = \rho_1, u = 1$,
coarsened, truncated 
$X_1 = (50, 50.001, 50.002, ..., 100)$: 
rows 1, 2: logarithmically scaled polygons plotted against
filter index 512, 128, 64, and 16.
row 3: polygon lengths for
filter indices $\lfloor 2^n \rfloor$ (left) and
$\lfloor 1.11^n \rfloor$ (right)
}
\label{seqfig10r1c1fi512.png}
\end{figure}

\begin{figure}[!htbp]
\centering
\includegraphics[scale=.45]{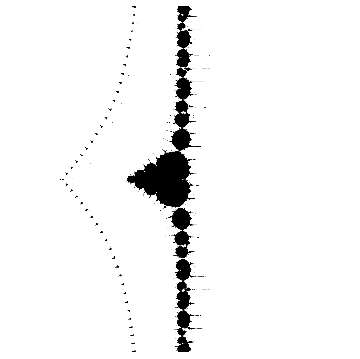}
\hskip .05in
\includegraphics[scale=.45]{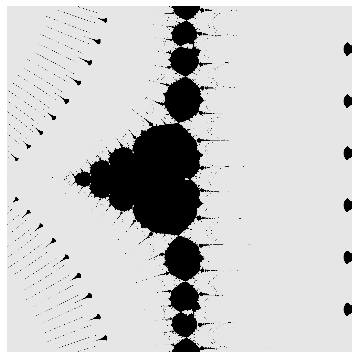}
\vskip .07in
\caption{left: $A_{\phi}$;
right: $\bf{C} - A_{\infty}$
}
\label{leaf1point1.png}
\end{figure}

For convenience, we reprint the conjecture:
\newline \newline
For each choice
of $\rho$,
there is a continuous function
$f^{cut}: \bf{C}^{\it cut} 
\rightarrow \bf{C}$
such that the restriction of 
$f^{cut}$ to $\vec{R}_u$
is a function $f$ with the
following properties:
\newline
(1) 
$f: \vec{R}_u \rightarrow \bf{C}$ 
is continuous and one-to-one,
\newline
(2) $f(z) \in \Phi_{z}$
for each $z \in \vec{R}_u$,
\newline
(3)  $f(0) = \psi_{\rho}$,
\newline
(4)  $\lim_{z \to \infty} 
f(z) =  \rho$,
\newline
(5) the image of $\vec{R}_u$ 
under $f$ is a nearly
logarithmic
spiral $S_{u,\rho}$ with center 
$\rho$,
\newline
(6)  among the sequences
$\vec{\phi} = 
(\phi_0, \phi_1, ...) \in 
\vec{\Phi}_{X,u}$
converging to $\rho$, 
all the points of one of them
(say, $\vec{\phi}_{u,\rho}$)
are interpolated by $S_{u,\rho}$;
the initial element of  
$\vec{\phi}_{u,\rho}$
is  
$\psi_{\rho}$.
\newline
(\it n.b. \rm Our notation suppresses the dependence of
$f$ on $u$ and $\rho$, and
the dependence of $f^{cut}$
on $\rho$.)
\newline \newline 
\textsection 
\thinspace
We are proposing this because 
a natural explanation
for our observation that
convergent
sequences $\vec{\phi} \in \Phi_{X, u}$
are apparently always nearly logarithmic
is that the tail of
such a $\vec{\phi}$ lies on 
a nearly logarithmic spiral.
The naturalness disappears 
(in our opinion) unless,
for a given $u$,
the spiral is independent of 
the choice of $X$.
Therefore, 
clauses (1) and (4)
both say something stronger,
namely, that
every single point of  $S_{u,\psi,\rho}$
is, in fact, an element of $\Phi_{x,u}$
for some $x \geq 0$.
(For the sake of clarity, we have
used
some redundancy in the 
statement of this conjecture.)
If we knew that such an $f$ exists,
at least, independent of the choice
of $X$'s with rational 
common differences, it seems 
clear that it would extend 
by continuity to the 
family of $X$'s with 
real common differences.
If $X_1, X_2$ with 
common differences $d_1, d_2$ both
rational are two instances of $X$
such that neither is a refinement
of the other, they clearly
have a refinement
in common with common difference a
rational number again. 
\newline \newline
Thus
(it seems to us) 
we can test conjecture 3 by
examining a tower of arithmetic
progressions $X_n$ such that 
$m > n \Rightarrow X_m$ is
a refinement of $X_n$.
(This is not the $X_n$ of 
section 1.2)
For a given $u$ each such $X$
determines (by way, for example, of the 
\it Mathematica 
\tt LinearModelFit \rm command)
a linear model of the data 
$(\theta, \log r)$ derived from
the points of a sequence 
$\vec{\phi} \in \Phi_{X,u}$
converging to a zero $\rho$
as follows: for $z_n \in \vec{\phi},
r = |z_n - \rho|$ and 
$\theta = \theta_{\vec{\phi}}(z_n)$
(see section 1.1.)
\newline \newline
In row 1 of Figure \ref{seqfig9r1c1.png}
we have plotted the slopes and intercepts of
such models for $u = 1$ and 
$\rho = \rho_1$
against an index on the horizontal
axis that specifies the refinement of a
particular $X = X_1 = (0, \frac 1{1000}, \frac 2{1000}, ..., 100)$
from $x = 50$ to $x = 100$.
(Actually,
for the sake of
efficiency, 
we worked in the
opposite direction: from
more highly-refined arithmetic
progressions to coarser ones.)
In row 2, we have plotted
the logarithms of the absolute
values of the
successive differences of the corresponding 
values from row 1. The common
differences of the $X$'s in the
tower increase from left to right; thus,
if the point corresponding to an
arithmetic progression $X_1$
lies to the left of the point
corresponding to an arithmetic progression
$X_2$, then $X_1$ is a refinement of $X_2$.
The row 2 values indicate
that the row 1 heights are a Cauchy sequence, 
thus a convergent sequence,
and we propose that the limits
of these two Cauchy sequences are the 
slope and intercept
of a log-linear model of $S_{u, \psi, \rho}$, 
\it i.e., \rm that $S_{u, \psi,\rho}$ is quite
possibly an exactly logarithmic spiral, but 
is, at least, very probably a nearly logarithmic 
spiral.
\newline \newline
As we have mentioned,
we did not, in fact, 
begin with a coarse
$X$ and refine it repeatedly 
to create a tower
of $X$'s.
That
procedure, to avoid redundant 
operations, would have 
required computing the
new fixed points at
each stage and then inserting
them into the sequence
determined from the
previous $X$. This seemed
too baroque.
Instead,
we began with the ``highly-refined'' 
arithmetic progression $X_1$ 
and chose from it successively less-refined
subsequences by the following process;
at stage $n$ of the process,
we selected every $n^{th}$ member
of the original 
$X_1$ to form a ``coarsening''
$X_n$ of $X_1$.  
Consequently, we were working
with subsequences $\vec{\phi}_n$ 
of the original 
$\vec{\phi}_1 = \vec{\phi} \in \Phi_{X_1,u}$,
one for each
of the coarsened arithmetic
progressions $X_n$.
The data we were going to model then
comprised (initially) 
a subsequence of the list of
pairs $(\theta, \log r)$ derived from 
$\vec{\phi}_1$. But,
because we were working with subsequences, 
the values of $\theta$ now violated
the minimality condition in clause (3)
of the definition of 
$\theta_{\vec{\phi}_n} = 
\theta_{(z_0, z_1, ...)}$ (say);
so we set $\theta_{\vec{\phi}_n}(z_0)= 
\arg (z_0 - \rho)$ as per clause (1),
and, for $k > 0$, chose 
$\theta_{\vec{\phi}_n}(z_k)$
according to clauses (2), (3) as applied,
not to the first $k$ elements of $\vec{\phi}_1$,
but to the first $k$ elements of $\vec{\phi}_n$.
Figure \ref{seqfig9r1c1.png} shows the plots
resulting from these operations
for $\rho = \rho_1, u = 1$,
and coarsenings $X_n$ of 
$X_1 = (0, \frac 1{1000}, \frac 2{1000}, ..., 100)$
from $x = 50$ to $x = 100$.
(In the caption to \ref{seqfig9r1c1.png},
we refer to this $n$ as the ``filter index.'')
We ignore the first $50,000$ members of
$X_1$ because the plots of $\log r$ vs. $\theta$
in this region are visibly non-linear
(even though, because only the tails count in 
this matter, the sequences in question do qualify as 
nearly logarithmic). We speculate
that the reason
for this observation is that 
\it Mathematica\rm's \tt FindRoot \rm command,
when searching for 
$\phi_{x,u} \in \Phi_{x,u}$ for relatively small $x$,
has not yet found the special sequence
mentioned in clause (5) of the conjecture--because,
for small $x, \phi_{x,u}$ is not yet close
enough to  
$\rho = \lim \vec{\phi}_1$. Nevertheless,
$\phi_{0,u}$, which by definition is a zeta fixed point,
is probably (we conjecture) so close to $\rho$ that it
is the zeta fixed point we have denoted as $\psi_{\rho}$;
this is the reason for clause (2) of conjecture 3.
\vskip .1in
\hskip -.23in
\thinspace
Figure \ref{seqfig10r1c1fi512.png}
shows the spirals (actually: polygons)
induced by the filter indices $512$ (row 1, column 1),
$128$ (row 1, column 2), $64$ (row 2, column 1),
and $16$ (row 2, column 2.) As this parameter decreases, the set of points included in
the plot grows and the polygons appear
to better approximate a spiral,
which we propose is the
$S_{1,\psi_{\rho},\rho}$ conjectured
to exist in conjecture 3. The
fourth row displays two plots of the
arc lengths of the polygons corresponding
to filter index 
$\lfloor 2^n \rfloor, 0 \leq n \leq 13$
(column 1) and to filter index
$\lfloor 1.1^n \rfloor, 0 \leq n \leq 100$
(column 2.) The arc lengths appear
to converge to what we propose is the arc length 
of the part of $S_{1,\psi_{\rho},\rho}$
with endpoints corresponding to 
$x = 50$ and $x = 100$. The convergence
(we wish to argue) is more evidence for
the existence of the spiral $S_{1,\psi_{\rho},\rho}$.
\section{\sc appendices}
\subsection{Extension of conjecture 3 to arbitrary
zeta fixed points.}
\begin{proposal}
For each choice
of $\rho$ and $\psi$
there is a continuous function
$f^{cut}: \bf{C}^{\it cut} 
\rightarrow \bf{C}$
such that the restriction of 
$f^{cut}$ to $\vec{R}_u$
is a function $f$ with the
following properties:
\newline
(1) 
$f: \vec{R}_u \rightarrow \bf{C}$ 
is continuous and one-to-one,
\newline
(2) $f(z) \in \Phi_{z}$
for each $z \in \vec{R}_u$,
\newline
(3)  $f(0) = \psi$,
\newline
(4)  $\lim_{z \to \infty} 
f(z) =  \rho$,
\newline
(5) the image of $\vec{R}_u$ 
under $f$ is a nearly
logarithmic
spiral $S_{u,\psi,\rho}$ with center 
$\rho$,
\newline
(6)  among the sequences
$\vec{\phi} = 
(\phi_0, \phi_1, ...) \in 
\vec{\Phi}_{X,u}$
converging to $\rho$, 
all the points of one of them
(say, $\vec{\phi}_{u,\psi,\rho}$)
are interpolated by $S_{u,\psi,\rho}$;
the initial element of  
$\vec{\phi}_{u,\psi,\rho}$
is  
$\psi$.
\newline
(\it n.b. \rm For readability,
our notation suppresses the dependence of
$f$ on $u, \psi$, and $\rho$, and
the dependence of $f^{cut}$
on $\psi$ and $\rho$.)
\end{proposal}
\subsection{Extension of conjecture 4 to arbitrary zeta fixed
points.}
\begin{proposal}
For each $X$ and $u$,
there is a function $G$ with the following
properties:\newline
(1) $G$ is meromorphic, independent 
of $\psi$,
and independent of $\rho$.
\newline
(2) The sequence  $\vec{\phi}_{u,\psi,\rho}$
from conjecture 3
satisfies
$G(\phi_n) =  \phi_{n-1}$
for  $n = 1, 2, ....$ 
\newline
(3)  Each $\rho$ is a repelling fixed point
of $G$,
\newline
(4) $G$ is many-to-one,
but for each $\rho$ and  each
$\psi$ there 
is a function 
$F_{\psi,\rho}$
such that 

(i)  $F_{\psi,\rho}^{(-1)} = G$,

(ii) consequently (in view of clause 2) 
for $n = 0, 1, ...$,
$$F_{\psi,\rho}(\phi_n) =  \phi_{n+1},$$

(iii) 
 $\rho$ is
an attracting fixed point of 
$F_{\psi,\rho}$,
and 

(iv) $F_{\psi,\rho}$ 
carries $S_{u,\psi,\rho}$
into itself.
\newline
(\it n.b. \rm Again for readability,
our notation suppresses the 
dependence of $G$,
and the dependence of
$F_{\psi, \rho}$ on $X$ and $u$. 
\end{proposal}
\newpage
\subsection{Two basins of attraction.}
\subsection{Figure 3.2 of \texorpdfstring{\cite{Br2}}{TEXT}.}
\begin{figure}[!htbp]
\centering
\includegraphics[scale=.45]{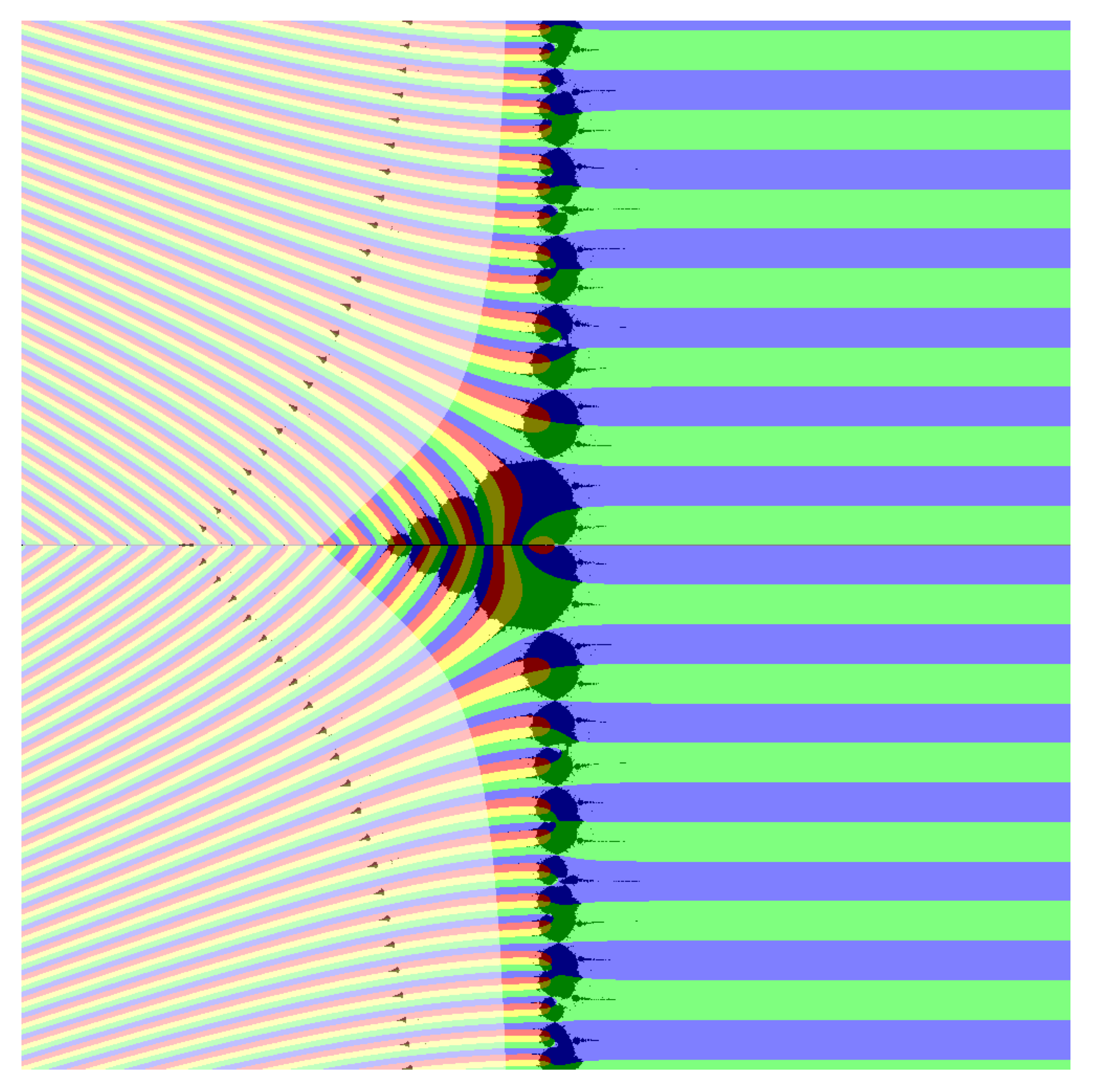}
\hskip .05in
\includegraphics[scale=.45]{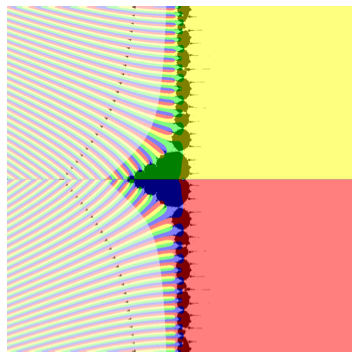}
\vskip .07in
\caption{left: zeros of zeta; right: zeros of
$s \mapsto \zeta(s) - s)$ (\it i.e.\rm, zeta
fixed points);
both superposed on a plot 
of $A_{\phi}$
}
\label{leaf3point2left.png}
\end{figure}

This figure is reproduced below as 
Figure \ref{leaf3point2left.png}.
\newpage
\bibliography{bibtexcite}
\href{mailto:barrybrent@member.ams.org}{barrybrent@member.ams.org}
\end{document}